\documentclass[reqno,12pt,twoside]{amsart}
\newtheorem{theorem} {Theorem}

\def\[{\left[}
\def\]{\right]}
\def\({\left(}
\def\){\right)}

\DeclareMathOperator{\lhs}{LHS}
\newcommand\pbin[2]{\genfrac{\langle}{\rangle}{0pt}{}{#1}{#2}}

\begin{document}
\author{Theresia Eisenk\"olbl}

\title{Proof of a partition identity conjectured by Lassalle}

\address{Institut f\"ur Mathematik der Universit\"at Wien\\
Strudlhofgasse 4, A-1090 Wien, Austria.}

\email{teisenko@radon.mat.univie.ac.at}

\begin{abstract}
We prove a partition identity conjectured by Lassalle
(Adv\@. in Appl\@. Math\@. {\bf 21} (1998), 457--472). 
\end{abstract}

\maketitle

The purpose of this note is to prove the theorem below which was
conjectured by Lassalle \cite{Lass1,Lass2}. In order to state the
theorem, we introduce the following notations.
Let $(a)_n = a(a+1)\cdots(a+n-1)$. For a partition $\mu$ of $n$ let
the length $l(\mu)$ be the number of the parts of $\mu$, $m_i$ the
number of parts $i$,
$z_{\mu}=\prod _{i\ge 1} ^{}{i^{m_i(\mu)}m_i(\mu)!}$ and
$\pbin{\mu}{r}$ the number of
ways to choose $r$ different cells from the diagram of the partition
$\mu$ taking at least one cell from each row. 
Then the following
theorem holds for $n \ge 1$.

\begin{theorem}
\begin{multline} \label{main}
\sum _{\vert \mu \vert =n} ^{}{\pbin{\mu}{r}\frac {X^{l(\mu)-1}}
{z_{\mu}}\sum _{i=1} ^{l(\mu)}{(\mu_i)_s}}\\
=(s-1)!\binom{n+s-1}{n-r}\[\binom{X+r+s-1}{r}-\binom{X+r-1}{r}\]
\end{multline}
\end{theorem}

\begin{proof}
We first observe that $\prod _{i\ge 1} ^{}{i^{m_i(\mu)}}=\prod _{i=1}
^{l(\mu)}{\mu_i}$ and that $\frac {l(\mu)!} {m_1!\cdots m_n!}$ is the
number of compositions of $n$ which are permutations of the parts of
$\mu$. Let us denote this number by $C(\mu)$.
After division by $s!$ the left-hand side can be rewritten as
\begin{align*}
\frac {\lhs}{s!}&=\sum _{\vert \mu \vert =n} ^{}{C(\mu)
\pbin{\mu}{r}\frac {X^{l(\mu)-1}}
{l(\mu)! \prod _{i=1} ^{l(\mu)}{\mu_i}}\sum _{i=1}
^{l(\mu)}{\binom{\mu_i+s-1}{s}}}\\
&=\sum_{l=1}^{\infty}{\sum _{\underset{\mu_j\ge
1}{\mu_1+\dots+\mu_l=n}} ^{}{\frac {X^{l-1}}
{l!\mu_1\cdots\mu_l}\pbin{\mu}{r}\sum _{i=1} ^{l}{\binom{\mu_i+s-1}{s}}}}
\end{align*}
For the composition $\mu$, $\pbin{\mu}{r}$ counts the ways of choosing
$r$ points in the diagram of the composition. If we choose $r_i$
points from part $\mu_i$, there are $\prod _{i=1}
^{l}{\binom{\mu_i}{r_i}}$ possible choices. 
Summing over all possible compositions $r=r_1+\dots +r_l$, where
every part is $\ge 1$ gives $\pbin{\mu}{r}$.
Thus we get for the left-hand side of \eqref{main}
\begin{multline*}
\frac {\lhs} {s!}
=\sum _{l=1} ^{\infty}{
\sum_{\underset{\mu_j\ge 1}{\mu_1+\dots+\mu_l=n}} 
    {\frac {X^{l-1}} {l!} \sum_{\underset{r_j\ge 1}{r_1+\dots +r_l=r}}
{\frac {1} {r_1\cdots r_l}
 \binom{\mu_1-1}{r_1-1} \cdots \binom{\mu_l-1}{r_l-1}
    \sum _{i=1}^{l}{\binom{\mu_i+s-1}{s}}
}
    }
                     }
\end{multline*}
It is easy to see that $\binom{\mu_i+s-1}{\mu_i-1}
\binom{\mu_i-1}{r_i-1}=(-1)^{r_i-1}\binom{-s-1}{r_i-1}\binom{\mu_i+s-1}
{r_i+s-1}$.
Now we can evaluate the sum over the $\mu_j$ by repeated application
of the Chu-Vandermonde summation formula:
$$\sum _{\mu_1+\dots+\mu_l=n} ^{}{\binom{\mu_1-1}{r_1-1}\cdots
\binom{\mu_l-1}{r_l-1}\binom{\mu_i+s-1}{s}}=
(-1)^{r_i-1}\binom{-s-1}{r_i-1}\binom{n+s-1}{r+s-1}.$$
Thus, we get for the left-hand side of \eqref{main}
\begin{equation} \label{c}
\frac {\lhs} {s!}=
\sum _{l=1} ^{\infty}{\frac {X^{l-1}} {l!}
     \sum_{\underset{r_j\ge 1}{r_1+\dots +r_l=r}}
{\frac {1} {r_1\cdots r_l}
 \sum _{i=1}^{l}{(-1)^{r_i-1}\binom{-s-1}{r_i-1}\binom{n+s-1}{r+s-1}}} }.
\end{equation}
The factor $\binom{n+s-1}{r+s-1}=\binom{n+s-1}{n-r}$ can be taken
outside of all the sums. By comparison of \eqref{main} and \eqref{c},
we see that it remains to prove

\begin{multline}\sum _{l=1} ^{\infty}{\frac {X^{l-1}} {l!}
     \sum_{\underset{r_j\ge 1}{r_1+\dots +r_l=r}}
{\frac {1} {r_1\cdots r_l}
 \sum _{i=1}^{l}{(-1)^{r_i-1}\binom{-s-1}{r_i-1}}}}\\
=
\frac {1} {s}\[\binom{X+r+s-1}{r}-\binom{X+r-1}{r}\].
\end{multline}
This can be done by using generating functions. We multiply both
sides of the
equation by $\Phi^r$ and sum over all $r\ge 0$. The right-hand side can
be evaluated by the binomial theorem and gives
\begin{equation} \label{rsgen}
\frac {1} {s}\( (1-\Phi)^{-X-s}-(1-\Phi)^{-X} \).
\end{equation}
For the left-hand side we need the power series expansion of the
logarithm and the equation
$$\sum _{r_i=1} ^{\infty}{\binom{r_i+s-1}{s}\frac {\Phi^{r_i}}
{r_i}}=\frac {1} {s} ((1-\Phi)^{-s}-1),
$$
which can be derived from the binomial theorem.
So the generating function corresponding to the left-hand side of
\eqref{rsgen} evaluates as follows:

\begin{multline*} 
\sum _{l=1} ^{\infty}{\frac {X^{l-1}} {l!}\sum _{r_1=1} ^{\infty}{\frac
{\Phi^{r_1}} {r_1}\sum _{r_2=1} ^{\infty}{\frac {\Phi^{r_2}}
{r_2}{\cdots\sum _{r_l=1} ^{\infty}{\frac {\Phi^{r_l}} {r_l}{\sum
_{i=1} ^{l}{\binom{r_i+s-1}{s}}}}}}}}\\
\begin{split}
&=\sum _{l=1} ^{\infty}{\frac {X^{l-1}} {l!}\sum _{i=1} ^{l}{\(\log \frac
{1} {1-\Phi}\)^{l-1}\frac {1} {s}\(\(1-\Phi\)^{-s}-1\)}}\\
&=\frac {1} {s}((1-\Phi)^{-s}-1)\sum _{l=1} ^{\infty}{\frac {\(X\log \frac
{1} {1-\Phi}\)^{l-1}} {(l-1)!}}\\
&=\frac {1} {s}((1-\Phi)^{-s}-1)e^{X\log\frac {1} {1-\Phi}}\\
&=\frac {1} {s}((1-\Phi)^{-s}-1)(1-\Phi)^{-X}\\
&=\frac {1} {s}((1-\Phi)^{-X-s}-(1-\Phi)^{-X}).
\end{split}
\end{multline*}

This is equal to \eqref{rsgen}, so the theorem is proved.
\end{proof}

\end{document}